\RequirePackage{fix-cm}
\documentclass[smallextended]{svjour3}       
\smartqed  
\usepackage{graphicx}
\usepackage{mathptmx}      
\usepackage{latexsym}
\usepackage{amssymb,amsmath,latexsym}
\usepackage{cite,natbib}
\usepackage{graphicx}
\usepackage{epsfig}
\usepackage{color}
\begin{document}

\title{Reaction-diffusion waves of blood coagulation
}


\author{Tatiana Galochkina         \and
        Anass Bouchnita \and
        Polina Kurbatova \and
        Vitaly Volpert 
}


\institute{T. Galochkina \at
              Camille Jordan Institute, University Lyon 1, Villeurbanne, 69622, France\\
Department of Biophysics, Faculty of Biology, M.V. Lomonosov Moscow State University, Leninskie gory 1, Moscow, 119992, Russia\\
Federal Research Clinical Center of Federal Medical \& Biological Agency of Russia, Orekhovy boulevard 28, Moscow, 115682, Russia \\
              Tel.: +7-495-9391116\\
              Fax: +7-495-939115\\
              \email{tat.galochkina@gmail.com}           
           \and
           A. Bouchnita \at
              Camille Jordan Institute, University Lyon 1, Villeurbanne, 69622, France\\
Laboratoire de Biom\'etrie et Biologie Evolutive, UMR 5558 CNRS, University Lyon 1, Lyon, 69376, France\\
Laboratory of Study and Research in Applied Mathematics, Mohammadia School of Engineers,  Mohamed V university, Rabat, Morocco
                      \and
           P. Kurbatova \at
Laboratoire de Biom\'etrie et Biologie Evolutive, UMR 5558 CNRS, University Lyon 1,
Lyon, 69376, France
                       \and
           V. Volpert \at
Camille Jordan Institute, University Lyon 1, Villeurbanne, 69622, France
}

\date{Received: date / Accepted: date}

\maketitle

\begin{abstract}
One of the main characteristics of blood coagulation is the speed of clot growth. This parameter strongly depends on the speed of propagation of the thrombin concentration in blood plasma. In the current work we consider mathematical model of the coagulation cascade and study the existence, stability and speed of propagation of the reaction-diffusion waves of blood coagulation. We also develop a simplified one-equation model that describes the main features of the thrombin wave propagation. For this equation we estimate the wave speed analytically. The resulting formulas give a good approximation for the speed of wave propagation in a more complex model as well as for the experimental data.

\end{abstract}

\section{Introduction}	

\subsection{Mathematical models of blood coagulation}

The main function of the coagulation system is terminating bleeding, caused by the injury of the vessel wall, with a fibrin clot.
The reaction of fibrin polymerization appears at the final stage of the proteolytic enzymatic cascade where the activated clotting factors act as catalyst for activation of the others \citep{Butenas2001,Orfeo2005}. Mature form of fibrin molecules can aggregate into the long branching fibers and form a complex network which serves as thrombus scaffold. The key enzyme of the coagulation cascade is thrombin as it catalyzes fibrinogen conversion to fibrin and its concentration has a crucial influence on the kinetics of the clot formation \citep{Hemker1993,Hemker1995,Butenas2001}.
To prevent the spontaneous formation of thrombi the activation reactions are normally regulated by the action of plasma inhibitors \citep{Pieters1988,Koppelman1995,Monkovic1990,Panteleev2002}. The balance between coagulation and anti-coagulation systems is important for the normal organism functioning and any alternations can lead to severe pathological states: thrombosis or, on the contrary, disseminative bleeding \citep{Colman2006,Askari2010}.

The influence of different factors on the coagulation process was studied both experimentally and using multiple theoretical approaches. As compared with the experiment, parameters in the theoretical studies can be varied much easier allowing to detect not only experimentally observed regimes of blood coagulation \citep{Stortelder1997,Leiderman2011,Krasotkina2000,Ataullakhanov1998,Bouchnita2015} but also to suppose their possible variations for the conditions that are hard to reproduce in the experiment \citep{Zarnitsina2001a}. Model results also provide data about the possible spatiotemporal distribution of all the blood factors participating in the coagulation cascade, while the main parameter used to measure the dynamics of the clot growth experimentally is fibrin clot density\citep{Ataullakhanov1998,Krasotkina2000,Ataullakhanov2002,Panteleev2006,Ovasenov2005}. Recently, the hypothesized in multiple mathematical models shape of the thrombin concentration profile in blood plasma was obtained in the direct experiment \citep{Dashkevich2012}.

The biochemical reactions in plasma are normally described with a set of ordinary or partial differential equations \citep{Stortelder1997,Leiderman2011,Krasotkina2000,Jones1994,Butenas2004,Lacroix2012} while the blood flow is modeled using continuous, discrete or mixed methods.  
Despite the differences of the modeling approaches, the systems of equations used for the description of the biochemical reactions are rather close, and they can usually be reduced to some simplified models \citep{Zarnitsina2001a,Ataullakhanov1998,Bouchnita2015}. In Section 2 of the current paper we consider a reaction-diffusion model of the intrinsic pathway of blood coagulation cascade and perform its reduction under certain assumptions. For the reduced system we prove a number of properties of its solutions (Section 3,4) that allow us to replace the analysis of the whole system by the analysis of one equation on the dynamics of thrombin concentration.

\begin{figure}[h]
\centerline{\includegraphics[scale=0.6]{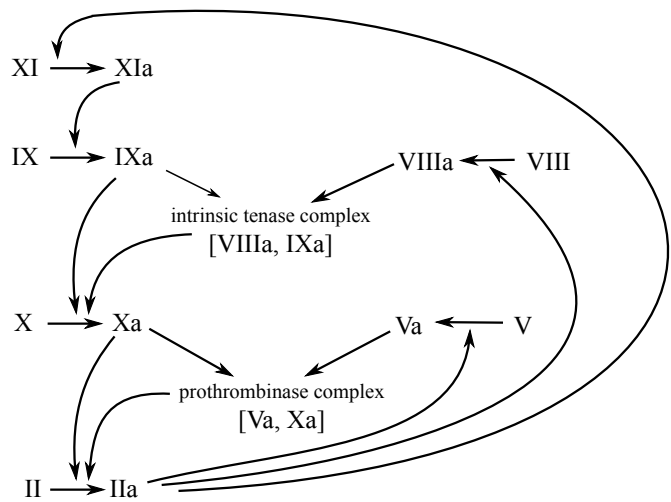}}
\caption{The main activation reactions of the intrinsic pathway of the coagulation cascade. \label{Scheme}}
\end{figure}

\subsection{The speed of thrombin propagation}

Thrombin is the key enzyme of the coagulation cascade as it catalyzes fibrinogen cleavage to fibrin which in turn forms haemostatic clot. Thrombin formation is a result of the prothrombin activation that can be initiated by different blood factors.
Normally the process is launched by the tissue factor that is expressed on the surfaces of all human cells except the endothelium and is bared in case of the vessel wall damage (extrinsic pathway), or through the activation of the factor XI after the contact with the foreign surface (contact activation) \citep{Orfeo2005,Orfeo2008,Gailani1991}.
Both pathways lead to the activation of factor X \citep{Orfeo2005} that contributes to the thrombin activation.
Once the thrombin concentration reaches the threshold value, further prothrombin activation takes place due to the positive feedback loops of the coagulation cascade (intrinsic pathway) \citep{Orfeo2008,Orfeo2005,Panteleev2006}. Thrombin controls activation of factor XI \citep{Gailani1991} and also of factors V \citep{Monkovic1990} and VIII whose activated forms (Va, VIIIa) increase activity of factors Xa and IXa by formation of the prothrombinase and intrinsic kinase complexes respectively \citep{Butenas2004,Orfeo2005,Baugh1996,Scandura1996} (Fig.~\ref{Scheme}).

As the result of these biochemical reactions, the thrombin wave propagates in the direction from the injury site to the vascular lumen. According to the experimental data, after thrombin reaches substantial concentration the speed of the clot growth does not anymore depend on the way of the initial activation of the coagulation system \citep{Orfeo2008,Ovasenov2005} and thrombin wave profile stays constant in time \citep{Butenas2004,Ataullakhanov1998,Tokarev2006}. The speed of propagation of the thrombin concentration serves then as one of the important criteria for the validation of the mathematical models of the coagulation system. Moreover, this speed depends strongly on the activity of the blood factors and can be significantly lowered or increased in case of illnesses disturbing the normal functioning of the coagulation system \citep{Colman2006,Askari2010}. This phenomena can be simulated using mathematical models of the whole coagulation cascade  \citep{Tokarev2006,Zarnitsina1996a,Lobanov2005} or the combination of analytical and numerical approaches as it was done by \cite{Pogorelova2014}.

In Section 5 of the current work we propose a new theoretical approach for the estimation of the speed of the thrombin wave propagation. After reducing the system of PDE describing the coagulation cascade to one equation, we build the analytical estimates on the speed of the propagation of the traveling wave solution using two different methods. We then compare the estimates given by analytical formulas with computational values of the speed as well as with the experimental data.

\section{Mathematical model}

\subsection{Complete model}

The main biochemical reactions of the coagulation cascade can be described by the following system of partial differential equations: 

\begin{equation}
\label{full}
\begin{aligned}
\frac{\partial U_{11}}{\partial t}  & =  D \Delta U_{11}  + r_{11} V_{11} T - h_{11}U_{11}, \\
\frac{\partial V_{11}}{\partial t}& = D \Delta V_{11} - r_{11} V_{11} T , \\
\frac{\partial T}{\partial t}& = D \Delta T  + r_{2} U_{10} \dfrac{P}{P + K_{2m}} + \overline{r_2}C_1\dfrac{P}{P + \overline{K_{2m}}} - h_2 T, \\
\frac{\partial P}{\partial t} &= D \Delta P  - r_{2} U_{10} \dfrac{P}{P + K_{2m}} - \overline{r_2}C_1\dfrac{P}{P + \overline{K_{2m}}}, \\
\frac{\partial C_1}{\partial t}& = D \Delta C_1 + k_{510}U_{5}U_{10} - h_{510}C_1,\\
\frac{\partial U_{5}}{\partial t} &= D \Delta U_{5} + r_{5} V_{5}T  - h_{5}U_{5},\\
\frac{\partial V_{5}}{\partial t} &= D \Delta V_{5}  - r_{5} V_{5} T,\\
\frac{\partial U_{10}}{\partial t} &= D \Delta U_{10}  + r_{10} V_{10}U_{9} + \overline{r_{10}}V_{10} C_2  - q_{10} U_{10},\\
\frac{\partial V_{10}}{\partial t} &=D \Delta V_{10}  - r_{10} V_{10}U_{9} - \overline{r_{10}} V_{10} C_2,\\
\frac{\partial C_2}{\partial t}& = D \Delta C_2 + k_{89}U_{8}U_{9} - h_{89}C_2, \\
\frac{\partial U_{8}}{\partial t} &= D \Delta U_{8}  + r_{8} V_{8}T - q_{8}U_{8}, \\
\frac{\partial V_{8}}{\partial t} &= D \Delta V_{8}  - r_{8} V_{8}T , \\
\frac{\partial U_{9}}{\partial t} &=D \Delta U_{9}  + r_{9} V_{9}U_{11} - q_{9} U_{9}, \\
\frac{\partial V_{9}}{\partial t} &= D \Delta V_{9}  - r_{9} V_{9}U_{11} .
\end{aligned}
\end{equation}
Here, $U_i$ and $V_i$ denote the concentrations of the activated and inactivated forms of the $i$-th factor respectively, $C_1$ denotes the concentration of the prothrombinase complex and $C_2$ corresponds to the concentration of the intrinsic tenase complex. The first term in each equation describes diffusion of the factor in blood plasma. The process of factor activation depends both on concentrations of its inactivated form and on the concentration of the activating factor, thus we consider the factor activation as the bimolecular reactions with the rate constants $r_i$. The inhibition of the activated factors is assumed to be the first order reaction with the rate constants $q_i$. Reaction of the thrombin activation obeys to the Michaelis-Menten mechanism with Xa and prothrombinase complex playing the role of enzyme and $K_{2m}, \; \overline{K_{2m}}$ being the Michaelis constants \citep{Zarnitsina2001a,Zarnitsina1996,Anand2006}. The reactions of complex formation are considered as the bimolecular reactions. The rate constants of the complex formation are $k_i$; as the rates of the complex degradation are much lower, we neglect the corresponding terms. In this work we do not consider the reaction of the fibrinogen conversion to fibrin as it does not influence the other reactions of the cascade and thus can be decoupled from the system. The resulting model consists of fourteen PDEs and it takes into account the main reactions of the intrinsic pathway~\eqref{full}.


\subsection{Model reduction}

\begin{figure}[h]
\centerline{
\includegraphics[scale=0.6]{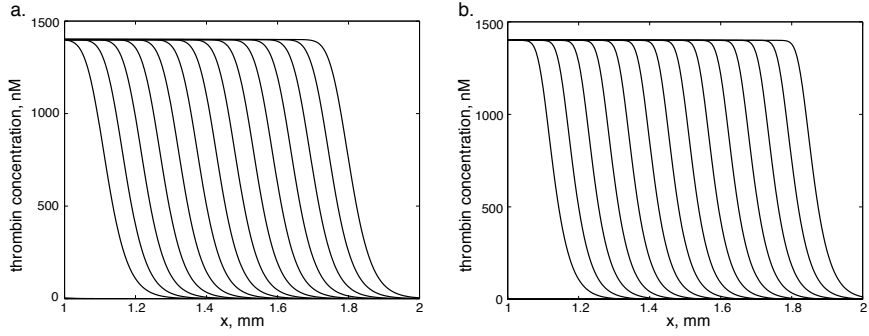}}
\caption{Propagation of thrombin wave for the model with prothrombin (left) and for the reduced model (\ref{reduced}) (right). Each next concentration profile is drawn after 1 min of physical time, the speed of the wave propagation is about $0.05$~mm/min. \label{Zarn}}
\end{figure}

To perform further analysis of system \eqref{full} we will reduce the number of equations.
First, let us apply the assumption of the detailed equilibrium to the complex formation equations and set:

\begin{equation}
C_1 = \dfrac{k_{510}}{h_{510}}U_{10}U_{5}, \; C_2 = \dfrac{k_{89}}{h_{89}}U_{9}U_{8}.
\end{equation}
This approximation is justified if the reaction constants $k_{510}, h_{510}, k_{89}, h_{89}$ are sufficiently large.
Next, we assume that the concentrations of the inactivated factors remain constant. This assumption implies that these factors are in excess and they do not limit the rate of the corresponding reactions. It allows us to consider the activation reactions as first-order reactions with $k_i = r_iV_i^0$. 
Then, we replace the Michaelis-Menten kinetics by the corresponding first-order reaction assuming that the constants $K_{2m}, \; \overline {K_{2m}}$ are
sufficiently large. Finally, we replace $P$ in the equation for $T$ using the approximation $P=P_0-T$. This relation is exact if $h_2=0$
and it gives a good approximation for small $h_2$ (Figure \ref{Zarn}).
The reduced system has the form:

\begin{equation}
\begin{aligned}\label{reduced}
\frac{\partial T}{\partial t} &= D \Delta T  + \left(k_{2} U_{10} + \overline{k_2}\dfrac{k_{510}}{h_{510}}U_{10}U_{5}\right)\left(1 - \dfrac{T}{T_0}\right) - h_2T, \\
\frac{\partial U_{5}}{\partial t} &= D \Delta U_{5}  + k_{5} T - h_{5}U_{5} ,\\
\frac{\partial U_{8}}{\partial t}& = D \Delta U_{8}  + k_{8} T - h_{8}U_{8},\\
\frac{\partial U_{9}}{\partial t}& = D \Delta U_{9}  + k_{9} U_{11} - h_{9}U_{9},\\
\frac{\partial U_{10}}{\partial t} &= D \Delta U_{10}  + k_{10} U_{9}  + \overline{k_{10}} \dfrac{k_{89}}{h_{89}}U_{9}U_{8}   - h_{10}U_{10},\\
\frac{\partial U_{11}}{\partial t} &= D \Delta U_{11}  + k_{11} T - h_{11}U_{11},
\end{aligned}
\end{equation}
where
$\overline{k_{10}} = r_{10}U_{10}^0, \; \overline{k_{2}} = r_2T_0$.

A similar model was previously considered by \cite{Zarnitsina1996a,Anand2006}. The parameter values of this model were verified using the experimental data. We take these values for the current study (Table~\ref{Param}, {Appendix C}). The thrombin wave propagates with a constant velocity of about $0.05$~mm/min (Fig.~\ref{Zarn}) that is in a good agreement with the experimental data \citep{Tokarev2006}.


\section{Existence and stability of waves}

Let us set $u = (T, U_5, U_8, U_9, U_{10}, U_{11})$. Then system (\ref{reduced}) can be written in the vector form:

\begin{equation}
\frac{\partial u}{\partial t} = D \Delta u  + F(u) , \label{tw}
\end{equation}
where $F=(F_1,...,F_6)$ is the vector of reaction rates in equations (\ref{reduced}).
It satisfies the following property:

$$ \dfrac{\partial F_i}{\partial u_j} \ge 0, \; \forall i \ne j . $$
This class of systems is called the monotone systems. These systems have a number of properties similar to those for one scalar equation including the maximum principle. It allows the proof of existence and stability of the wave solutions for such systems and estimation of the speed of the wave propagation \citep{VolpertEPDE}. 
In order to apply these results we need to begin with the analysis of the existence and stability of the stationary points of system~\eqref{reduced}.

\subsection{Stationary points of the kinetic system}

Consider the system of ordinary differential equations:

\begin{equation}\label{ode}
\dfrac{du}{dt} = F(u)
\end{equation} 
Its equilibrium points satisfy the following relations:

\begin{gather}
\label{g1}
U_{5} = \dfrac{k_{5}}{h_{5}}T, \; U_{8} = \dfrac{k_{8}}{h_{8}}T, \; U_{11} = \dfrac{k_{11}}{h_{11}}T, \; U_9 = \dfrac{k_9k_{11}}{h_9h_{11}}T, \\ 
U_{10} = \dfrac{k_9k_{11}}{h_{10}h_9h_{11}}\left(k_{10}T + \overline{k_{10}}\dfrac{k_{89}}{h_{89}}T^2\right) ,
\label{g2}
\end{gather}
where $T$ is a solution of the equation $P(T)=0$. Here $P(T)=aT^4 + bT^3 + cT^2 + dT$,

\begin{gather*}
a = \dfrac{\overline{k_{10}}k_{89}k_8\overline{k_2}k_5k_{510}k_9k_{11}}{h_{89}h_{8}h_5h_{10}h_{510}h_{9}h_{11}}, \;\;\;\;
d = -\dfrac{k_2k_{10}k_9k_{11}}{h_9h_{11}h_{10}} + h_2T_0 , \\
b = -\dfrac{\overline{k_{10}}k_{89}k_8\overline{k_2}k_5k_{510}k_9k_{11}}{h_{89}h_{8}h_5h_{10}h_{510}h_{9}h_{11}}T_0 + \dfrac{k_{10}\overline{k_2}k_5k_{510}k_9k_{11}}{h_5h_{10}h_{510}h_{9}h_{11}} + \dfrac{\overline{k_2k_{10}}k_{89}k_8k_9k_{11}}{h_{89}h_{8}h_{9}h_{11}}, \\
c = -\dfrac{k_{10}\overline{k_2}k_5k_{510}k_9k_{11}}{h_5h_{10}h_{510}h_{9}h_{11}}T_0 + \dfrac{k_2k_{10}k_{9}k_{11}}{h_9h_{11}} - \dfrac{k_2k_{89}k_8k_9k_{11}}{h_{89}h_8h_9h_{11}h_{10}}T_0 .
\end{gather*}
Hence stationary points of system \eqref{ode} can be found through the stationary points $T^*$  of the equation

\begin{equation}
\label{eq1d}
\dfrac{\partial T}{\partial t} = - P(T)
\end{equation}
and equalities (\ref{g1}), (\ref{g2}).

Let us determine the number of positive roots of the polynomial $P(T)$. We set
$P(T) = TQ(T)$, where $Q(T) = aT^3 + bT^2 + cT + d$. 
The number of positive roots of $Q(T)$ can be found as follows. First, we consider a function $Q'(T) = 3aT^2 + 2bT + c$. If it has no zeros, then $Q(T)$ is increasing and has one positive root if and only if $Q(0) < 0$.
Otherwise, we denote by $T_1, T_2$ the nonzero solutions of the equation  $Q'(T) = 0$:
$ T_{1,2} = (-b \pm \sqrt{b^2 - 3ac})/(3a)$.
Then the polynomial $Q(u)$ has one positive root in one of the cases:
\begin{itemize}
\item $T_1 \leq 0, \; Q(0) < 0$,
\item $0 \leq T_1 < T_2, \; Q(0) < 0 \; {\rm and}  \; Q(T_1) > 0, \; Q(T_2) > 0 \; {\rm or} \; Q(T_1) < 0$
\end{itemize}
and it has two positive roots if $0 < T_2, \; Q(0) > 0, \; Q(T_2) < 0$.

Stability of stationary points of the system \eqref{ode} can be determined from the stability of stationary points of equation (\ref{eq1d}).
The following theorem holds (see {Appendix~A} for the proof).

\begin{theorem}
There is one to one correspondence between stationary solutions $u^* = (T^*, U_5^*, $ $U_8^*, U_9^*, U_{10}^*, U_{11}^*)$
of system \eqref{reduced} and $T^*$ of equation \eqref{eq1d} given
by (\ref{g1}), (\ref{g2}). The principal eigenvalue of the matrix $F'(u^*)$ is positive (negative) if and only if $P'(T^*) < 0$ ($P'(T^*) > 0$).
\end{theorem}

Thus, we can make the following conclusions about the existence and stability of stationary points of the kinetic system of equation~\eqref{ode}. It always has the trivial solution $u^*=0$. It has one (two) positive solution if and only if the polynomial $P(T)$ has one (two) positive roots. A positive solution $u^*$ is stable if and only if $P'(T^*)>0$.


\subsection{Wave existence and stability}

We can now formulate the theorem on the existence of waves for the system~\eqref{reduced}.

\begin{theorem}
Suppose that $P(T^*)=0$ for some $T^*>0$ and $P'(0) \neq 0,\; P'(T^*) \neq 0$. Let $u^* = (T^*, U_5^*, $ $U_8^*, U_9^*, U_{10}^*, U_{11}^*)$
be the corresponding stationary solutions of system (\ref{ode}) determined by relations (\ref{g1}), (\ref{g2}).
\begin{itemize}
\item Monostable case. If there are no other positive roots of the polynomial $P(T)$, then system (\ref{reduced}) has monotonically decreasing travelling wave solutions $u(x,t)=w(x-ct)$ with the limits $u(+\infty)=0, u(-\infty)=u^*$ for all values of the speed $c$ greater than or equal to the minimal
speed $c_0$,
\item Bistable case. If there is one more positive root of the polynomial $P(T)$ in the interval $0 < T < T^*$, then system (\ref{reduced}) has a monotonically decreasing travelling wave solutions $u(x,t)=w(x-ct)$ with the limits $u(+\infty)=0, u(-\infty)=u^*$ for a unique value of $c$.
\end{itemize}
\end{theorem}

The proof of the theorem follows from the general results on the existence of waves for monotone systems of equation \citep{VolpertEPDE,Volpert2000}.
Let us note that conditions on the stability of stationary points follow from the assumption of the Theorem~2 and Theorem~1. We have $P'(T^*)>0$ in both cases since it is the largest root of the polynomial increasing at infinity. The sign of $P'(0)$ is negative if there is no other root in between and
positive if there is one more root. 

Monotone travelling wave solutions of monotone systems are asymptotically stable \citep{VolpertEPDE,Volpert2000} that gives global stability in the bistable case. 
In the monostable case the wave is globally stable for the minimal speed $c_0$ and stable 
with respect to small perturbations in a weighted norm for $c>c_0$ \citep{Volpert2000}.

The unique wave speed in the bistable case and the minimal wave speed in the monostable case admit minimax representations. Below we will use them for the bistable system as this case is more appropriate for the applications considered in the current work. Indeed, travelling wave solution of system~\eqref{reduced} describes propagation of the thrombin concentration in the blood flow due to the reactions of the coagulation cascade. In this system the convergence to the travelling wave solution takes place only if the initial concentrations of blood factors exceed some critical level, otherwise the clot formation does not start because of the action of plasma inhibitors. Described dependency on the initial conditions and stability of zero solution correspond to the bistable case. In the monostable case, on the contrary, a small perturbation would result in the solution converging to the propagating wave. In terms of the coagulation system functioning, monostable case corresponds to the spontaneous disseminated coagulation blocking blood circulation.

Finally, let us note that in the Theorem~2 we consider only the case of a single positive root of the polynomial and the case of two positive roots. If $P(T)$ has three positive roots the system would be monostable with a stable intermediate stationary point. While this case is interesting from the point of view of wave existence and stability, it is less relevant for the modelling of blood coagulation, and we will not discuss it here.


\section{Speed of wave propagation}

One of the main objectives of this work is to obtain an analytical approximation of the wave speed for the blood
coagulation model \eqref{reduced}. We proceed in two steps. 
First, we reduce the system \eqref{reduced} to a single equation and justify this reduction. Then in the next section, we obtain some estimates of the wave speed for one reaction-diffusion equation.


\subsection{System reduction}

In order to simplify the presentation, we describe the method of reduction for the system of two equations:

\begin{align}
u'' + cu' + f(u,v) & = 0,\label{uu}\\
v'' + cv' + \dfrac{1}{\varepsilon}(au - bv) &= 0, \label{vv}
\end{align}
where $\varepsilon$ is a small parameter, $\dfrac{\partial f}{\partial v} > 0$ and the system \eqref{uu},\eqref{vv} is bistable.
If we multiply the second equation by $\varepsilon$ and take a formal limit as $\varepsilon \to 0$, then we have $v = \dfrac{a}{b}u$, and the first equation can be rewritten as follows:

\begin{equation}
u'' + cu' + f\left(u,\dfrac{a}{b}u\right) = 0 \label{tw-1}.
\end{equation}
Let us recall that the value of the speed $c=c_\varepsilon$ in system~(\ref{uu}), (\ref{vv}) and $c=c_0$ for the scalar equation (\ref{tw-1}) are unknown,
and in general they are different from each other. We will demonstrate that $c_\varepsilon \to c_0$ as $\varepsilon \to 0$:

\begin{theorem} The speed of wave propagation for system \eqref{uu},~\eqref{vv} converges to the speed of the wave propagation for equation~\eqref{tw-1} as $\varepsilon \to 0$.
\end{theorem}

Singular perturbations of travelling waves are extensively studied by \cite{VolpertEPDE}. Here we present another method of proof based on the estimates of the wave speed. This method is simpler and gives not only the limiting value of the speed for $\varepsilon = 0$ but also the estimates of the speed value for any positive $\varepsilon$. In the following section we describe the approach in details and construct the wave speed estimates for the system (\ref{uu}), (\ref{vv}).


\subsection{Wave speed estimate}

We get the following estimates from the minimax representation of the wave speed in the bistable case \citep{Volpert2000} :

\begin{equation}\label{mm}
  \min \left( \inf_x S_1(\rho), \inf_x S_2(\rho) \right) \leq c \leq \max \left( \sup_x S_1(\rho), \sup_x S_2(\rho) \right) ,
\end{equation}
where

$$ S_1(\rho) = \dfrac{\rho_1'' + f(\rho_1,\rho_2)}{-\rho_1'} \;, \;\;\;\; 
S_2(\rho) = \dfrac{\rho_2'' + (a\rho_1 - b\rho_2)/\varepsilon}{-\rho_2'} \; ,  $$
$\rho=(\rho_1, \rho_2)$ is an arbitrary test function continuous together with its second derivatives, monotonically decreasing (component-wise) and having the same limits at infinity as the wave, $\rho(+\infty)=0, \; \rho(-\infty)=u^*$.

Let us choose the following test functions:

\begin{gather}
\rho_1 = u_0, \;\;\;\;
\rho_2 = \dfrac{a}{b}u_0 - \varepsilon f\left(u_0, \dfrac{a}{b}u_0\right)\dfrac{a}{b^2} \; ,
\end{gather}
where $u_0$ is the solution of the equation \eqref{tw-1}. Neglecting the second-order terms with respect to $\varepsilon$, we get:

\begin{multline}
S_1(\rho) = 
\left(u_0'' + f\left(u_0, \dfrac{a}{b}u_0 - \varepsilon \dfrac{a}{b^2} f\left(u_0, \dfrac{a}{b}u_0\right)\right)\right)/(-u_0') = \\
\left( u_0'' + f\left(u_0, \dfrac{a}{b}u_0\right) - \varepsilon \dfrac{a}{b^2}f_v\left(u_0, \dfrac{a}{b}u_0\right)f\left(u_0, \dfrac{a}{b}u_0\right)\right) /(-u_0') =
 c_0 + \varepsilon \varphi(x), \label{test_1}
\end{multline}
where 

$$\varphi(x) = \dfrac{a}{b^2u_0'} f_v\left(u_0, \dfrac{a}{b}u_0\right)f\left(u_0, \dfrac{a}{b}u_0\right),$$
and $c_0$ is the value of the speed in \eqref{tw-1}. Next, 

\begin{equation}
\label{test_2}
 S_2(\rho) = \dfrac{u_0'' + f\left(u_0, \dfrac{a}{b}u_0\right) -
 \dfrac{\varepsilon}{b} \left(f\left(u_0, \dfrac{a}{b}u_0\right) \right)'' }{-u_0' + \dfrac{\varepsilon}{b} \left(f\left(u_0, \dfrac{a}{b}u_0\right)\right)'} \; = c_0 + \varepsilon \psi(x) ,
\end{equation}
where

$$ \psi =  \dfrac{c_0}{b u_0'} \left(f\left(u_0, \dfrac{a}{b}u_0\right)\right)' + \dfrac{1}{b u_0'} \left(f\left(u_0, \dfrac{a}{b}u_0\right)\right)'' . $$
Hence from \eqref{test_1}, \eqref{test_2} we obtain the estimate

\begin{equation}
c_0 + \varepsilon \max \left\{\min_{x}\varphi, \min_{x} \psi \right\} \le c \le c_0 + \varepsilon \min \left\{\max_{x}\varphi, \max_{x} \psi \right\},
\end{equation}
where $c_0$ is the wave propagation speed for \eqref{tw-1}, the functions $\varphi(x), \; \psi(x)$ are bounded. The proof of Theorem~3 follows from this estimate.


\section{One equation model}

\subsection{Reduction to the single equation}

\begin{figure}[h]
\begin{center}
\includegraphics[scale=0.6]{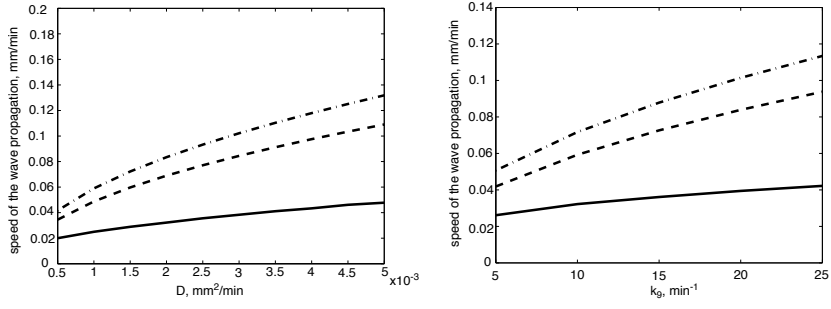}
\end{center}
\caption{Speed of wave propagation (mm/min) as a function of $D$ (left) and $k_9$ (right). Solid line: reduced model \eqref{reduced}; dashed line: two-equation model \eqref{2eq}; dash-dot line: one equation model \eqref{1eq}.\label{D}}
\end{figure}

If the reaction rate constants in the equations  of system (\ref{reduced}) for the variables
$U_{9}$, $U_{10}$, $U_{5}$ and $U_{8}$ are sufficiently large, then we can replace these equations by the following algebraic relations
(cf. Section 4.1):

\begin{gather*}
U_{5} = \dfrac{k_5}{h_5}T, \; U_{8} = \dfrac{k_8}{h_8}T, \; U_{9} = \dfrac{k_9}{h_9}U_{11}, \;
U_{10} = U_{11}\dfrac{k_9}{h_9h_{10}}\left(k_{10} + \dfrac{\overline{k_{10}}k_{89}}{h_{89}}\dfrac{k_8}{h_8}T\right).
\end{gather*}
Then instead of system (\ref{reduced}) we obtain the following system of two equations:

\begin{equation}
\begin{aligned} \label{2eq}
\frac{\partial T}{\partial t} &= D \Delta T  + U_{11}\dfrac{k_9}{h_9h_{10}}\left(k_{10} + \dfrac{\overline{k_{10}}k_{89}}{h_{89}}\dfrac{k_8}{h_8}T\right)\left(k_2 + \dfrac{\overline{k_2}k_{510}}{h_{510}}\dfrac{k_5}{h_5}T \right)\left(1 - \dfrac{T}{T_0}\right) - h_2T,\\
\frac{\partial U_{11}}{\partial t} &= D \Delta U_{11}  + k_{11} T - h_{11}U_{11}.
\end{aligned}
\end{equation}
Similarly, we can reduce this system to the single equation:

\begin{equation}\label{1eq}
\frac{\partial T}{\partial t} = D \Delta T + \dfrac{k_9k_{11}}{h_9h_{10}h_{11}}T\left(k_{10} + \dfrac{\overline{k_{10}}k_{89}}{h_{89}}\dfrac{k_8}{h_8}T \right)\left(k_2 + \dfrac{\overline{k_2}k_{510}}{h_{510}}\dfrac{k_5}{h_5}T\right) \left(1 - \dfrac{T}{T_{0}} \right) - h_2T.
\end{equation}
We realize this reduction in two steps in order to compare the one-equation model with the system~\eqref{reduced} as well as the intermediate model of two equations~\eqref{2eq}. Numerical simulations show that for the values of
parameters in the physiological range, all three models give the wave speed of the same order of magnitude (Fig.~\ref{D}).
The two equation model (\ref{2eq}) gives a better approximation of the model (\ref{reduced}) than the single equation (\ref{1eq}).
However, the latter demonstrates the same parameter dependence of the wave speed as other models. Taking into account the complexity of the initial model~\eqref{reduced}, the approximation provided by one equation is acceptable. In what follows we obtain analytical formulas for the wave speed for the one equation model.

\subsection{Dimensionless model}

In dimensionless variables

\begin{equation}
T = T_0u, \; t = \dfrac{\tilde{t}}{h_2}, \; D = \tilde{D}h_2.
\end{equation}
we write equation \eqref{1eq} in the following form:

\begin{equation}\label{nond}
\frac{\partial u}{\partial \tilde{t}} = \tilde{D} \Delta u + M_1u\left(1 + M_2u\right)\left(1 + M_3u\right)\left(1 - u\right) - u,
\end{equation}
where:
\begin{gather*}
M_1 = \dfrac{k_2k_9k_{10}k_{11}}{h_2h_9h_{10}}, \; M_2 = \dfrac{k_8k_{89}\overline{k_{10}}}{k_{10}h_8h_{89}}T_{0}, \; M_3 = \dfrac{\overline{k_2}k_5k_{510}}{k_2h_5h_{510}}T_{0}.
\end{gather*}
Analysis of the rate constant values allows us to further simplify the equation. As $M_3\gg1$ we can approximate equation~\eqref{nond} by the following equation:

\begin{equation}
\frac{\partial u}{\partial \tilde{t}} = \tilde{D} \Delta u + M_1M_3u^2\left(1 + M_2u\right)\left(1 - u\right) - u.
\end{equation}
Let us note that in the term $u_1^2(1 + M_2u_1)$ the first component corresponds to the prothrombin activation by the factor Xa and the second one corresponds to the prothrombin activation by the [Va, Xa] complex. Since during the propagation phase the rate of activation by prothrombinase complex is several orders of magnitude higher than the activation by Xa itself, we can neglect the first component. Thus, applying the assumption of the detailed equilibrium for the second equation, we finally obtain the following equation for the thrombin concentration:

\begin{equation}
\frac{\partial u_1}{\partial \tilde{t}} = \tilde{D} \Delta u_1 + b u_1^3\left(1 - u_1\right) - u_1, \label{u1}
\end{equation}

where:
\begin{equation}\label{bbb}
b = M_1M_2M_3.
\end{equation}


\subsection{Wave speed estimate}

The equation \eqref{u1} can be rewritten in the more general form:
\begin{equation}
\frac{\partial u}{\partial t} = D \Delta u + b u^n\left(1 - u\right) - \sigma u . \label{u}
\end{equation}
Travelling wave solution of~\eqref{u} satisfies the equation:

\begin{equation}
\label{w}
D w'' + c w' + b w^n(1-w) - \sigma w = 0 .
\end{equation}
Here we will present two approximate analytical methods to determine the wave speed.

\subsubsection{Narrow reaction zone method}

One of the methods to estimate the wave speed for the reaction-diffusion equation is the narrow reaction zone method developed in combustion theory \citep{Zeldovich1938}.
Let us rewrite the equation \eqref{w} in the form:

\begin{equation}
\label{narrow}
D w'' + c w' + F(w) - \sigma w = 0, \; F(w) = w^n(1-w)
\end{equation}
We assume that the reaction takes place at one point $x=0$ in the coordinates of the moving front. Then outside of the reaction zone we consider the linear equations:

\begin{equation}
\label{wout2}
\left\{
\begin{array}{lr}
D w'' + c w'  - \sigma w = 0, & x > 0, \\
D w'' + c w' = 0, & x <  0.
\end{array}
\right.
\end{equation}
These equations should be completed with the jump conditions at the reaction zone. In order to derive them, we omit the first derivative $w'$ at the reaction zone since it is small in comparison with two other terms:

\begin{equation}
D w'' + F(w) = 0. \label{react}
\end{equation}
Multiplying \eqref{react} by $w'$ and integrating through the reaction zone we obtain the following jump conditions:

\begin{equation}
\label{jump2}
(w'(+0))^2 - (w'(-0))^2  = \frac2D \int\limits_0^{w^*} F(w) dw ,
\end{equation}
considered together with the condition of the continuity of solution $w(+0) = w(-0)$. 

Solving \eqref{wout2} we have:

\begin{equation}
\label{wsol2}
w = \left\{ \begin{array}{lr} w_*, & x < 0, \\ w_*\exp\left(\dfrac{-c-\sqrt{c^2 + 4D\sigma}}{2D}\right), & x> 0. \end{array}\right.
\end{equation}
Then from \eqref{jump2} and \eqref{wsol2} we obtain the following equation for the wave speed:

\begin{equation}\label{eqcnar}
\hat{c}^2 + \hat{c}\sqrt{\hat{c}^2 + 4D\sigma} + 2D\sigma = A, \;\;\; A = \dfrac{4D}{w_*^2} \int\limits_0^{w^*} F(w) dw.
\end{equation}
Hence

\begin{equation}
\label{c2}
\hat{c}_1 = \dfrac{A - 2D\sigma}{\sqrt{2A}} \; , \;\;\;\;
A = 4bD\left(\dfrac{w_*^{n-1}}{n+1} - \dfrac{w_*^{n}}{n+2}\right).
\end{equation}

\begin{figure}[h]
\centerline{\includegraphics[scale=0.4]{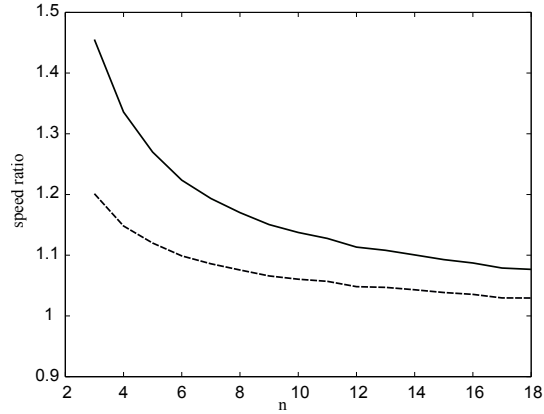}}
 \caption{Ratio of wave speeds found numerically and analytically for different values of $n$; $\sigma = 0.01, \; D = 2, \; b = 10$. Solid line: $\dfrac{c}{\hat{c}_1}$, dashed line $\dfrac{c}{\hat{c}_2}$.}
 \label{sp}
\end{figure}

This formula gives a good approximation of the wave speed found numerically for $n \geq 3$ (Fig.~\ref{sp}). The approximation improves with growth of $n$. The obtained formula gives an estimation of the speed from below (see {Appendix~B} for the justification of the method).

\subsubsection{Piecewise linear approximation}

Consider equation (\ref{narrow}) written in the form

$$D w'' + c w' + f(w) = 0 , $$
where $f(w) = w^n(1-w) - \sigma w$ and $f(0)=f(w_*)=0$.
Let us introduce the following approximation of this equation:
\begin{equation}
\label{pw}
D w'' + c w' + f_0(w) = 0,
\end{equation}
with
\begin{equation}
f_0(w) = \left\{ \begin{array}{lr} \alpha w, & 0<w<w_0, \\  \beta(w - w_*), & w_0<w<w_*,      \end{array}\right.
\end{equation}
where
\begin{equation}
\alpha = f'(0) , \; \beta = f'(w_*).
\end{equation}
In case of equation \eqref{u} we have:
\begin{equation}
\begin{aligned}
\alpha = -\sigma, \;\;\;\;
\beta = bnw_*^{n-1} - b(n+1)w_*^{n} - \sigma.
\end{aligned}
\end{equation}
We find the value of $w_0$ from the additional condition:

\begin{equation}
\int_0^{w_*} f(w)dw = \int_0^{w_*} f_0(w) dw .
\end{equation}
Hence we obtain the following equation with respect to $w_0$:
\begin{equation}\label{w0}
\dfrac{\alpha - \beta}{2}w_0^2 + \beta w_*w_0 + r = 0,
\end{equation}
where
\begin{equation}
r = -\beta w_*^2 - \int_0^{w_*} f(w)dw.
\end{equation}
Taking into account the explicit form of the function $f(w)$, we obtain:
\begin{equation}
r = bw_*^{n+1}\left(-\dfrac{n}{2} - \dfrac{b}{n+1}\right) + bw_*^{n+2}\left(\dfrac{n+1}{2} + \dfrac{1}{n+2} \right) + \sigma w_*^2.
\end{equation}
From \eqref{w0} we get:
\begin{equation}
w_0 = \dfrac{-\beta w_* + \sqrt{\beta^2w_*^2 - 2(\alpha - \beta)r}}{\alpha - \beta}.
\end{equation}

Thus, instead of \eqref{pw} we consider the following equation:
\begin{equation}
\left\{
\begin{array}{lr}  Dw'' + cw' + \beta (w - w_*) = 0, & x<0 ,\\ D w'' + cw' + \alpha w = 0, & x > 0, \end{array}
\right.
\end{equation}
with the additional conditions on the continuity of solution and its first derivative:
$$ w(0) = w_0, \;\;\; w'(-0) = w'(+0) . $$
We find the explicit solution:
\begin{equation}
\left\{
\begin{array}{lr} w = (w_0 - w_*)\exp\left(x\dfrac{\sqrt{c^2 - 4\beta D} - c}{2D}\right) + w_*, & x<0 ,\\
w = w_0\exp\left(x\dfrac{-\sqrt{c^2 - 4\alpha D} - c}{2D}\right) , & x > 0. \end{array}
\right.
\end{equation}
From the condition of continuity of the derivative we obtain the following formula:
\begin{equation}\label{chat2}
\hat{c}_2 = \dfrac{\sqrt{D}(\alpha \bar{w}^2 - \beta)}{\sqrt{(\bar{w} - 1) (\alpha \bar{w}^2 - \beta \bar{w})}}, \;\;\;\; 
\bar{w} = \dfrac{w_0}{w_0 - w_*} \; .
\end{equation}
It gives a good approximation of the wave speed for equation~\eqref{narrow}  (Fig.~\ref{sp}).

\subsection{Comparison of the estimated speed of the wave propagation with the complete model and experimental data}

\subsubsection{Comparison of the estimated speed with the computational speed in system~\eqref{reduced}}

Considering the system \eqref{reduced} and taking the parameter values for \eqref{c2},~\eqref{chat2} according to~\eqref{bbb}, we approximate the speed of wave propagation by the following formula obtained by the narrow reaction zone method:

\begin{equation}
\hat{c_1} = \sqrt{D}\dfrac{bT_0^2 - \dfrac45{bT_0^3} - 2h_2}{\sqrt{2\left(bT_0^2 - \dfrac45{bT_0^3}\right)}},
\end{equation}
where
$$
b = \dfrac{k_9k_{11}\overline{k_{10}}k_8k_{89}\overline{k_{2}}k_5k_{510}T_0^2}{h_9h_{10}h_{11}h_8h_{89}h_5h_{510}},
$$
and by the piecewise linear approximation:

\begin{equation}
\hat{c_2} = \dfrac{\sqrt{D}\left( - 3b{T_0}^2 -h_2\overline{T} + 4bT_0^3 - h_2\right)}{\sqrt{\left(T_0-1\right)\overline{T}\left(-h_2\overline{T} - 3bT_0^2 + 4bT_0^3 + h_2\right)}},
\end{equation}
where:
\begin{multline}
\hspace*{3cm}  \overline{T} = \dfrac{T_*}{T_* - T_0}, \;\;\;\;
T_* = \frac{-3bT_0^2 + 4bT_0^4 + h_2}{4bT_0^2 - 3bT_0} + \\
 \dfrac{\sqrt{\left(3bT_0^2 - 4bT_0^3 - h_2\right)^2 - 2b(4T_0 - 3)T_0^2\left(-\frac32bT_0^2 - \frac{b^2}{4}T_0^2 + \frac{11}{5}bT_0^3 + h_2\right)}}{4bT_0^2 - 3bT_0}.
\end{multline}

\begin{figure}[h]
\begin{center}
\includegraphics[scale=0.3]{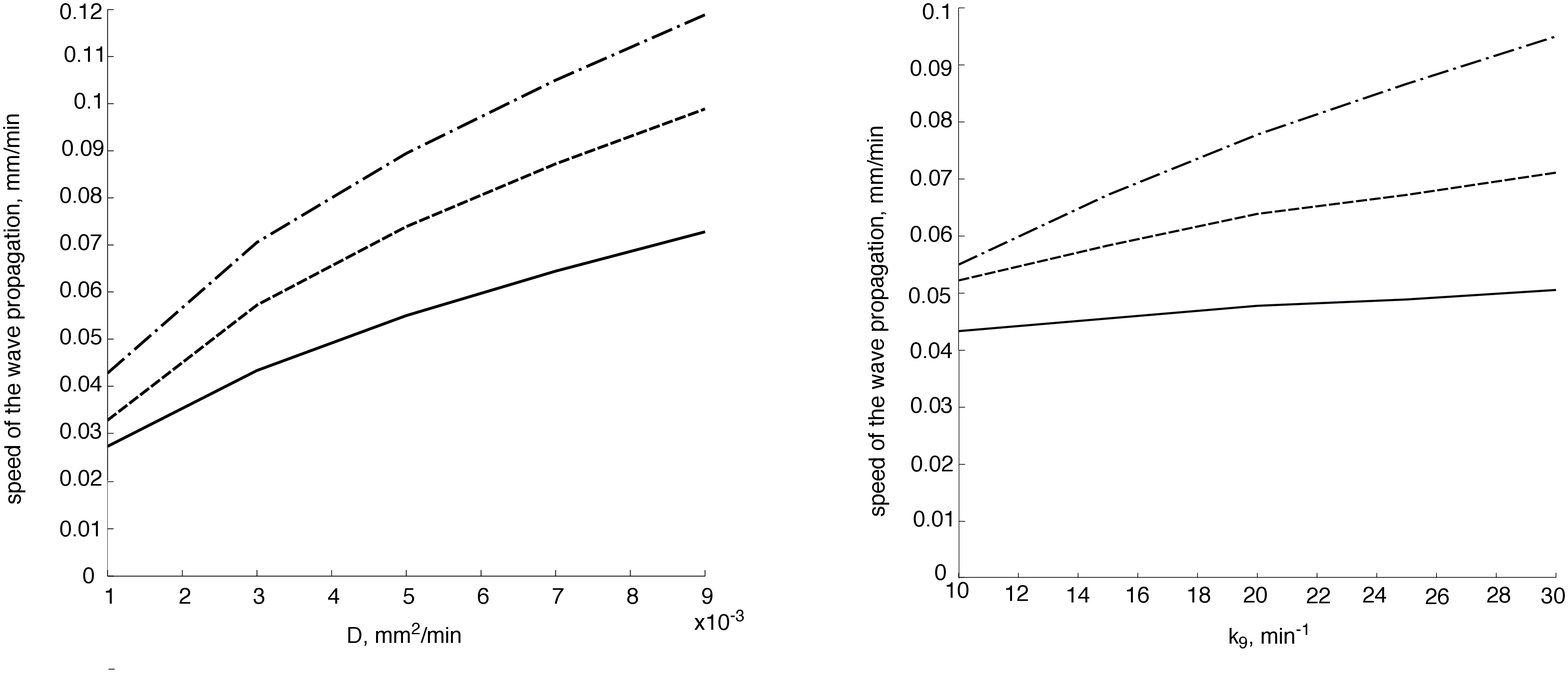}
\end{center}
\caption{Speeds of wave propagation (mm/min) as function of $D$ (left) and $k_9$ (right). Solid line:  model \eqref{reduced}; dashed line: narrow reaction zone approximation; dash-dot line: piecewise linear approximation.\label{Compl}}
\end{figure}

We compare the speed of wave propagation for model~\eqref{reduced} found numerically with the analytical formulas (Fig.~\ref{Compl}). As it was demonstrated above, the computational speed for the one-equation model is higher than for the complete model (Fig.~\ref{D}). The analytical formulas for the speed of the wave propagation for one-equation model give in turn the estimates from below (Fig.~\ref{sp}). As the result, the analytical estimates for one equation give a better approximation of the speed in the complete model than the numerical speed for one equation (Fig.~\ref{Compl}). If we then compare two different analytical estimates for the wave speed in one-equation model, we can conclude that narrow reaction zone method gives the speed further from the one-equation computational speed than piecewise linear approximation (Fig.~\ref{sp}) but at the same time it better approximates the wave speed in the complete model (the narrow reaction zone speed is $1.5$ times higher than the computational one).

\subsubsection{Comparison with experimental data}

\begin{figure}[h]
\centerline{\includegraphics[scale=0.6]{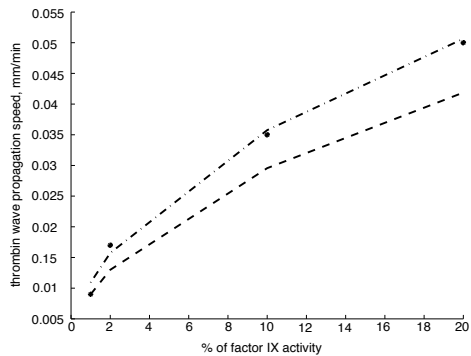}}
\caption{Speeds of the thrombin wave propagation (mm/min) as function of percentage of factor IX activity. Dots: experimental data; dashed line: narrow reaction zone approximation; dash-dot line: piecewise linear approximation. \label{experiment}}
\end{figure}

Analytical formulas for the wave propagation speed are in a good correlation with the experimental data. As an example, by \cite{Tokarev2006} we can find the results of the \emph{in vitro} experiments, measuring the speed of clot growth for the patients with hemophilia. This illness affects the activity of the factor IX thus slowing down the process of prothrombin conversion into thrombin. To compare the formula for the speed of wave propagation with the experimental data we fitted the parameters only for one experimental point corresponding to 1\% of factor IX activity. Then, we considered different values of the activity of factor IX which corresponds to the parameter $k_9$ in our estimates (Fig.~\ref{experiment}). Both formulas provide a good approximation of the experimental results.

\section{Discussion and conclusions}

Coagulation cascade is a complex network of chemical reactions whose functioning depends strongly on various factors. Numerous mathematical models were developed to describe the influence of different parameters on the details of clotting process. The key stage of the coagulation cascade defining the dynamics of the clot formation is cumulative thrombin production due to the intrinsic pathway functioning. In the most of mathematical models this process is described in terms of PDE on the concentrations of the involved factors and complexes. 
The traveling wave solutions of such systems correspond to the experimentally observed thrombin wave propagating from the vessel wall to the vascular lumen \citep{Dashkevich2012}. 

In the current work we proved existence and stability of traveling wave solutions for the system describing intrinsic pathway of blood coagulation cascade. The main property of the considered system that allowed us to study the traveling waves analytically is its monotonicity. Not all of the previously considered models satisfy this property, but most of them can be reduced to the monotone systems under some conditions. Moreover, the system reduction can be performed in different ways depending on the assumptions on the kinetics of coagulation cascade. In the upcoming work we will consider another method of the system reduction and will prove the existence and stability of wave solutions for the resulting system.

The most important parameter defining the dynamics of clot growth is the speed of the thrombin wave propagation, or in terms of the mathematical model, the speed of propagation of the reaction-diffusion wave. In the current work we obtain analytical formula for the speed of wave propagation in the model of blood coagulation. We reduce the system of equations to one equation and then determine the wave speed for this equation. The method of reduction is based on the minimax representation of the wave speed applicable for monotone reaction-diffusion systems.
Though one equation model gives the value of the speed higher than for the system of equations, the order of magnitude and the dependence on parameters remain the same.  Furthermore, the resulting estimates give a good approximation of the available experimental data.

Despite the general character of the methods used in this work, the developed approaches imply some limitations. First, in our model we considered only a part of the coagulation cascade (intrinsic pathway) without taking into account neither the initial activation, nor the role of the activated protein C. This simplification allowed us to perform more detailed analysis of the resulting model. The independence of the final speed of the thrombin wave propagation on the nature of the stimuli that launched the clotting process was previously studied in multiple works \citep{Orfeo2005,Orfeo2008,Tokarev2006}. The role of protein C can be also omitted if we consider the reaction on the distance from the vessel wall. Indeed, the process of the protein C activation occurs due to its reaction on the surface of the vessel wall after a substantial amount of thrombin is already formed \citep{Anand2008} and thus does not influence the speed of the thrombin propagation in the vascular lumen.

The next important assumption that makes possible model reduction is that certain constants of the kinetic reactions are sufficiently large. However, for the realistic values of the reaction rates this property is not always satisfied thus leading to the difference in the solutions for the complete and reduced models. We also observe this effect in our model, but as the resulting parameter dependence remains the same and the solutions are still rather close we take the reduced models as acceptable approximations of the complete one.

One more important limitation of the considered approach is the application of the narrow reaction zone method for analytical estimation of the wave propagation speed for one equation model. As this method was originally developed for the description of the flame front propagation in the combustion theory, the function describing the reaction should be close to the exponential function. In our work the reaction function is described with the polynomial of the third degree that makes the obtained estimate less precise. 

Finally, despite all the assumptions made the obtained analytical estimates give good approximation of both computational and experimental speed of the thrombin propagation. The system reduction, analysis of the existence and stability of waves and estimates of wave speed developed in this work can be further expanded on other models.  


\begin{acknowledgements}
The authors thank the French-Russian scholarship program ``Metchnikov'' for the opportunity of the international collaboration.
\end{acknowledgements}

\bibliographystyle{spbasic}      
\bibliography{WaveBib}   


\section*{Appendix}
\appendix

\section{Proof of the Theorem 1 \label{A}}

\begin{proof}

Along with the system system 

\begin{equation}\label{aa1}
  \frac{du}{dt} = F(u)
\end{equation}
consider the system

\begin{equation}\label{aa2}
  \frac{du}{dt} = F_\tau(u)
\end{equation}
which depends on the parameter $\tau \in [0,1]$. They differ only by the equation for $T$ which is considered now in the following form:

\begin{equation*}
 \frac{dT}{dt} = \left(\tau U_{10} + (1-\tau)\varphi_{10}(T)\right)\left(k_{2} + \overline{k_2}\dfrac{k_{510}}{h_{510}}\left(\tau U_{5} + (1 - \tau)\varphi_{5}(T)\right)\right)\left(1 - \dfrac{T}{T_0}\right) - h_2T.
\end{equation*}
Here the functions $\varphi_i(T)$ are determined by the equalities:

\begin{gather*}
\varphi_{11}(T) = \dfrac{k_{11}}{h_{11}}T, \;
\varphi_9(T) = \dfrac{k_9k_{11}}{h_9h_{11}}T, \;
\varphi_{5}(T) = \dfrac{k_{5}}{h_{5}}T, \;
\varphi_{8}(T) = \dfrac{k_{8}}{h_{8}}T,\\
\varphi_{10}(T) = \dfrac{k_9k_{11}}{h_{10}h_9h_{11}}\left(k_{10}T + \overline{k_{10}}\dfrac{k_{89}}{h_{89}}T^2\right).\\
\end{gather*}
We can express $U_i, \; i = 5,8,9,10,11$ as functions of $T$ from the corresponding equations in (\ref{aa1}) or, the same,
from (\ref{aa2}): $U_i = \varphi_i(T)$. Therefore the
solutions of the system of equations $F_{\tau}(T) = 0$ coincide with the solutions of the system $F(T) = 0$.

Thus, systems (\ref{aa1}) and (\ref{aa2}) have the same stationary solutions for all $\tau \in [0,1]$.
For $\tau=1$ these two systems coincide. For $\tau=0$ the equation for $T$ in (\ref{aa2}) does not depend on other variables.
This will allow us to determined the eigenvalues of the corresponding linearized matrix.

It can be verified by the direct calculations that det $F_\tau'(u^*)=0$ if and only if det $F'(u^*)=0$ for all $\tau \in [0,1]$.
Suppose that the latter is different from zero. Then the principal eigenvalue of the matrix $F_\tau'$, which is real and simple,
cannot change sign when $\tau$ changes from $0$ to $1$. Hence the sign of the principal eigenvalue of the matrix $F'(u^*)$
is the same as for the matrix $F_0'(u^*)$. This matrix has the form:

\begin{equation*}
F_0'(u^*) = 
\bordermatrix{ & T & U_5 & U_8 & U_{11} & U_9 & U_{10} \cr
T & - P'(T^*) & 0 & 0 & 0 & 0 & 0 \cr
U_5 & k_5 & -h_5 & 0 & 0 & 0 & 0\cr
U_8 & k_8 & 0 & -h_8 & 0 & 0 & 0 \cr
U_{11} & k_{11} & 0 & 0 & -h_{11} & 0 & 0  \cr
U_9 & 0 & 0 & 0 & k_9 & -h_9 & 0  \cr
U_{10} & 0 & 0 & \overline{k_{10}}\dfrac{k_{89}}{h_{89}}U_9^* & 0 & k_{10} + \overline{k_{10}}\dfrac{k_{89}}{h_{89}} U_8^*& -h_{10}}
\end{equation*}

The principal eigenvalue of this matrix is positive if $P'(T^*) < 0$ and negative if this inequality is opposite.

\end{proof}

\section{Justification of the narrow reaction zone method}

\begin{figure}[h]
\centerline{\includegraphics[scale=0.5]{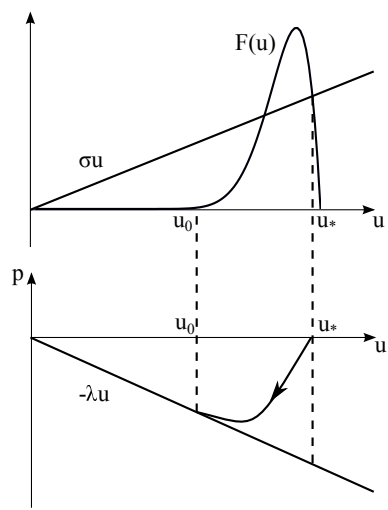}}
 \caption{Illustration of the narrow reaction zone method approximation.}
  \label{zf}
\end{figure}

Consider equation \eqref{narrow} and suppose for simplicity that $F(u)=0$ for $u \leq u_0$ and
$F(u)>0$ for $u_0 < u < 1$. Let $u^*$ be the maximal solution of the equation
$F(u)=\sigma u$ (Figure \ref{zf}). We will look for a decreasing solution of equation \eqref{narrow} with
the limits:

$$ u(-\infty)=u^*, \;\;\; u(+\infty)=0 . $$
Multiplying the equation \eqref{narrow} by $u'$ and integrate through the hole axis we obtain:

\begin{equation}\label{z2}
c = \dfrac{\int\limits_0^{u^*} F(u) du - \frac12 \sigma (u^*)^2}{\int\limits_{-\infty}^\infty (u'(x))^2 dx} \;.
\end{equation}
Along with equation \eqref{narrow} we consider the system of two first-order equations:

\begin{equation}\label{z3}
\left\{
\begin{array}{l} u'=p,\\ p' = - c p - F(u) + \sigma u .\end{array}
\right.
\end{equation}
The wave solution of \eqref{narrow} corresponds to the trajectory connecting the stationary points $(u^*,0)$ and $(0,0)$
(Figure \ref{zf}). This trajectory coincides with the line $p=\lambda u$ for $0 < u \leq u_0$,
where $\lambda$ is a negative solution of the equation

\begin{equation*}
\lambda^2 + c \lambda - \sigma = 0 .
\end{equation*}
The integral in the denominator of \eqref{z2} can be approximated by replacing the trajectory function by the straight line $p = -\lambda u$:

$$ \int_{-\infty}^\infty (u'(x))^2 dx = \int_0^{u^*} p(u) du \approx \frac12 \lambda (u^*)^2 . $$

Substituting this expression into \eqref{z2} we obtain the same formula for the speed as by the narrow reaction zone method \eqref{c2}.

Thus, narrow reaction zone method is equivalent to replacing the equation trajectory by the straight line.
Hence we can conclude that this method provides the estimate of the speed from below, and it also gives
asymptotically correct result in the limiting case as the support of the function $F(u)$ converges to a point.

\newpage
\section{Parameter values used for the simulations}

\begin{table}[ht]
\caption{Parameter rates used for the modeling of the coagulation cascade. \label{Param}}
{\begin{tabular}{@{}ccc@{}} 
parameter & value & units \\
\hline
$k_{11}$ & $0.000011$ & min$^{-1}$\\
$h_{11}$ & $0.5$ & min$^{-1}$\\
$k_{10}$ & $0.00033$ & min$^{-1}$\\
$\overline{k_{10}}$ & $500$ & min$^{-1}$\\
$h_{10}$ & $1$ & min$^{-1}$\\
$k_9$ & $20$ & min$^{-1}$\\
$h_9$ & $0.2$ & min$^{-1}$\\
$k_{89}$ & $100$ & nM$^{-1}$min$^{-1}$\\
$h_{89}$ & $100$ & min$^{-1}$\\
$k_8$ & $0.00001$ & min$^{-1}$\\
$h_8$ & $0.31$ & min$^{-1}$\\
$k_5$ & $0.17$ & min$^{-1}$\\
$h_5$ & $0.31$ & min$^{-1}$\\
$k_{510}$ & $100$ & nM$^{-1}$min$^{-1}$\\
$h_{510}$ & $100$ & min$^{-1}$\\
$k_2$ & $2.45$ & min$^{-1}$\\
$h_2$ & $2.3$ & min$^{-1}$\\
$K_{2m}$ & $58$ & nM\\
$\overline{K_{2m}}$ & $210$ & nM\\
$D$ & $0.0037$ & mm$^2$min$^{-1}$\\
$T_0$ & $1400$ & nM 
\end{tabular}}
\end{table}


%

\end{document}